\documentclass[12pt,reqno]{amsart}

\theoremstyle{fact}
\newtheorem{fact}{Fact}
\theoremstyle{plain}
\newtheorem{proposition}{Proposition}
\newtheorem{lemma}{Lemma}
\newtheorem{theorem}{Theorem}

\theoremstyle{definition}
\newtheorem{definition}{Definition}

\newcommand{\V}{\Vert}
\let\Sum=\sum\let\Lim=\lim\let\Sup=\sup\let\Inf=\inf\let\Max=\max
\def\lim{\Lim\limits}\def\sum{\Sum\limits}\def\inf{\Inf\limits}\def\sup{\Sup\limits}\def\max{\Max\limits}
\def\sss{\scriptscriptstyle}        
\def\ialign{\everycr{}\tabskip0pt plus0pt minus0pt\halign} 
\def\cases#1{\left\{\,\vcenter{\normalbaselines\mathsurround0pt\ialign{$##\hfil$&\quad##\hfil\crcr#1\crcr}}\right.}
\def\seq#1#2#3{\{#1_#2\}_{#2=#3}^\infty}
\def\sseq#1#2#3#4{\{#1_{#2_#3}\}_{#3=#4}^\infty}
\def\a{\alpha}\def\X{\frak X}\def\Y{\frak Y}\def\V{\Vert}
\def\I{\mathop{\rm I_{\raise1pt\hbox{$\sss1\kern.5pt2$}}}}
\def\In{\mathop{\rm I^n_{\raise1pt\hbox{\fiverm1\kern.5pt2}}}}  
\def\II#1{\mathop{\rm I_{\raise1pt\hbox{$\sss1\kern.5pt \sss #1$}}}}
\def\IIN{\mathop{\rm I_{\raise1pt\hbox{\fiverm2\kern.5pt2}}}}
\def\HI{\mathop{\rm\hat I_{\raise1pt\hbox{$\sss1\kern.5pt2$}}}}
\def\ts{\textstyle}
\begin{document}
\baselineskip=18pt
\title{Two characterizations of the standard unit vector basis of $l_1$}
\author{David\ Mitra}
\thanks{This work is a subset of the results contained in the author's dissertation, {\it Sequences that are unconditionally basic in 
        both $l_1$ and $l_2$}, which was written under the direction of Professor Maria Girardi.} 
\address{Department of Mathematics, University of South Carolina\vskip0pt Columbia, SC 29208, U.S.A.}
\email{mitra@math.sc.edu}
\date{September 19, 2000}
\begin{abstract}
      We show that for a sequence in a Banach space, the property of being stable under large perturbations characterizes the property of 
      being equivalent to the unit vector basis of $l_1$. We show that a normalized unconditional basic sequence in $l_1$ that is 
      semi-normalized in $l_\infty$ is equivalent to the standard unit vector basis of~$l_1$.
\end{abstract}
\maketitle
\section{Introduction}
In this note, we present two results concerning the standard unit vector basis of~$l_1$. The first result, Theorem~1, characterizes, up
to equivalence,  the  standard unit vector basis of $l_1$ as the unique sequence satisfying a certain stability property. This result is a generalization of
\cite[Theorem~II.11.5]{S}. Our second result, Theorem~2,  characterizes, among the sequences in~$l_1$ that are semi-normalized in both $l_1$ and 
$l_\infty$, those that are unconditionally basic
as simply those   that are equivalent to the standard unit vector basis of $l_1$. In contrast to Theorem~2, we 
present Proposition~1, which shows that not every copy of the standard unit vector basis of $l_1$ in $l_1$ is semi-normalized in $l_\infty$.

Throughout this note, we shall use the following notations: $\X$ and $\Y$ denote infinite dimensional, real Banach spaces; the symbol 
$\II p$ denotes the formal identity operator from $l_1$ into $l_p$; we denote the standard unit vector basis of $l_1$ by $\seq{e^1}n1$; 
finally, we denote the $i^{\rm th}$-coordinate of a sequence of real numbers $x$ by $x(i)$. All other notation and terminology, not 
otherwise explained, are as in~\cite{LT}.

We now recall the definitions pertinent to the foregoing discussion.
\begin{definition} 
The {\it unconditional basis constant\/} of an unconditional basic sequence $\seq xn1$ in $\X$ is the smallest constant~$K_1$ 
satisfying, for any sequence of signs $\{\epsilon_n\}_{n=1}^m$ and for any sequence of scalars $\{\a_n\}_{n=1}^m$, the inequality 
$$
  {\biggl\V\sum_{n=1}^m\epsilon_n\a_nx_n\biggr\V_\X\le K_1\biggl\V\sum_{n=1}^m\a_nx_n\biggr\V_\X}.
$$
\end{definition}

\begin{definition}  
A sequence $\seq xn1$ in $\X$ is {\it$(K_1,K_2)$-equivalent to the sequence\/} $\seq yn1$ in~$\X$, provided that
there exists a pair of positive numbers~$K_1,K_2$ such that
$$
  K_1\biggl\V\sum_{n=1}^m\a_n y_n\biggr\V_\X\le\biggl\V\sum_{n=1}^m\a_nx_n\biggr\V_\X\le K_2\biggl\V\sum_{n=1}^m\a_ny_n\biggr\V_\X,
$$
for each sequence $\{\a_n\}_{n=1}^m$ of scalars. We say that $\seq xn1$ is {\it equivalent} to $\seq yn1$ if it is
$(K_1,K_2)$-equivalent to $\seq yn1$ for some pair of positive numbers $K_1$, $K_2$. 
\end{definition}
\section{The Results}
Our first result concerns perturbations of basic sequences equivalent to $\seq {e^1}n1$. We first recall the classical result, due to 
Bessaga and Pe\l czy\'nski (see~\cite{BP}), on perturbations of basic sequences, which states that a basic sequence is ``stable under small perturbations''.

\begin{fact} 
Let $\seq xn1$ be a basic sequence in~$\X$ with coefficient functionals $\seq{x^*}n1$. If $\seq yn1$ is a sequence in $\X$ 
satisfying 
$$
  \sum_{n=1}^\infty\V x^*_n\V_{[x_i]^*_\X}\V x_n-y_n\V_{\X}<1,
$$ 
then $\seq yn1$ is a basic sequence, equivalent to~$\seq xn1$.
\end{fact}

In \cite[Theorem II.11.5]{S}, it is proved that $\seq{e^1}n1$ is actually ``stable under large perturbations''. To clarify this vague 
statement, and to facilitate our discussion, we make the following definition.

\begin{definition}
A sequence $\seq xn1$ in a Banach space $\X$ is a {\it $\delta$-\break dominated perturbation\/} of a sequence $\seq yn1$ in $\X$ provided
\begin{equation}\label{eq: wwww}
      \sup_{n\in\mathbb N}\V x_n-y_n\V_\X<\delta.
\end{equation}
\end{definition}

The result in \cite{S} mentioned above states that any $1$-dominated perturbation of $\seq {e^1}n1$  is equivalent to 
$\seq {e^1}n1$. According to Singer, this result is contained in the paper \cite{M} of V. D. Milman, who, in turn, attributes the 
result to V. I. Gurari\/\u\i. In Theorem~1, we show that, in fact, Definition~3 may be used to characterize those
sequences that are equivalent to $\seq{e^1}n1$.

\begin{theorem}\label{rm: bigpert} 
Let $\seq xn1$ be a normalized sequence in $\X$. The following conditions on $\seq xn1$ are equivalent. 
\begin{enumerate}
\item[{\rm 1)}] $\seq xn1$ is equivalent to the standard unit vector basis of $l_1$. 
\item[{\rm 2)}] There exists a positive number $\delta$ such that $\seq xn1$ is uniformly equivalent to any $\delta$-dominated perturbation of itself\/;
               that is, such that the constants $K_1$ and $K_2$ appearing in Definition~2 depend
               only on the supremum appearing in equation~{\rm(1)}.
\item[{\rm 3)}] There exists a positive number $\delta$ such that $\seq xn1$is equivalent to any $\delta$-dominated perturbation of itself.
\item[{\rm 4)}] $l_1$ is isomorphically embedded in $\X$, and there exists a positive number $\delta$ such that $\seq xn1$is equivalent to any 
               $\delta$-dominated perturbation of itself.
\end{enumerate}
\end{theorem}

\begin{proof} To see that condition~1) implies condition~2),
suppose (as we may) that $\seq xn1$ is $(k,1)$-equivalent to $\seq {e^1}n1$. Take $\delta=k$ and let $\seq yn1$ be a $\delta$-dominated 
perturbation of $\seq xn1$. Then $\seq yn1$ is norm bounded; and so, since $\seq xn1$ is equivalent to $\seq {e^1}n1$, we may define a bounded linear 
operator $T\colon[x_n]_\X\rightarrow[y_n]_\X$ by setting
$$
  T\biggl(\sum_{n=1}^\infty\a_n x_n\biggr)=\sum_{n=1}^\infty\a_n y_n,\qquad{\rm for\ }x=\sum_{n=1}^\infty\a_nx_n.
$$ 
Moreover, for $x=\sum_{n=1}^\infty \alpha_n x_n$, we have, setting $m=\sup_{n\in\mathbb N}\V x_n-y_n\V_\X$, that
$$
  \V(I-T)x\V_\X=\biggl\V\sum_{n=1}^\infty\a_n(x_n-y_n)\biggr\V_\X 
  \le\sup_{n\in\mathbb N}\V x_n-y_n\V_\X\sum_{n=1}^\infty|\a_n|\le{m\over k}\V x\V_\X.
$$
Since $mk^{-1}< 1$, by a well-known result, $T$ is an isomorphism from $[x_n]_\X$ onto $[y_n]_\X$ satisfying
$$
   {k\over k+m}\biggl\V\sum_{n=1}^\infty\a_n y_n\biggr\V_\X\le\biggl\V\sum_{i=1}^\infty 
   \a_nx_n\biggr\V_\X\le {k\over k-m}\biggl\V\sum_{n=1}^\infty\a_n y_n\biggr\V_\X,     
$$
for each $x=\sum_{n=1}^\infty\alpha_nx_n$, as desired. 

It is obvious that condition 2) implies condition 3).

To see that condition 3) implies the first (and non-trivial) part of condition 4), first assume that $\seq xn1$ has a weakly Cauchy subsequence 
$\sseq xnk1$ and let $\seq zn1$ be an arbitrary normalized sequence in $\X$. Let $\seq yn1$ be the $\delta$-dominated perturbation  of 
$\seq xn1$ defined by
$$
  y_j= \cases{x_{n_k}-{\delta\over2} z_k,&if $j=n_k$;\cr x_j,&otherwise.\cr}
$$
Then, by hypothesis, $\seq yn1$ is equivalent to $\seq xn1$.  From this, it follows that the linear map 
$\widetilde T\colon\{x_n\}_{n=1}^\infty\rightarrow\{y_n\}_{i=1}^\infty$ given by $\widetilde Tx_n=y_n$ for each $n\in\mathbb N$ extends to a well-defined linear 
isomorphism $T\colon[x_n]_\X\rightarrow[y_n]_\Y$.
In particular, the map $T$ is weak-weak continuous. Thus, since $\sseq xnk1$ is weakly Cauchy, so is
$\sseq ynk1=\{x_{n_k}-\delta z_k\}_{k=1}^\infty$. But then we arrive at the contradiction that the arbitrary normalized sequence 
$\{z_n\}_{n=1}^\infty$ is weakly Cauchy. It must be that $\seq xn1$ has no weakly Cauchy subsequence, whence the result follows by 
Rosenthal's $l_1$-theorem \cite{R}.

Towards showing that condition 4) implies condition 1), let $\seq zn1$ be a normalized sequence in $\X$ equivalent to $\seq {e^1}n1$. 
Let $\seq yn1$ be the $\delta$-dominated perturbation  of $\seq xn1$ defined by
$$
  y_n=x_n+ \ts{\delta\over2}z_n.
$$
Then, by hypothesis, $\seq xn1$ is $(k_1,k_2)$-equivalent to $\seq yn1$ for some pair of positive numbers $k_1,k_2$. Thus, for any 
sequence $\{\alpha_n\}_{n=1}^m$ of scalars we have
$$
  \biggl\V\sum_{n=1}^m\a_n x_n\biggr\V_\X\ge k_1 \biggl\V\sum_{n=1}^m \a_n y_n\biggr\V_\X
  \ge k_1\Biggl(\biggl\V\sum_{n=1}^m \ts{\delta\over2}\,\a_nz_n\biggr\V_\X-\biggl\V\sum_{n=1}^m\a_nx_n\biggr\V_\X\Biggr);
$$
whence,
$$
  {2(1+k_1)\over \delta\cdot k_1}\biggl\V\sum_{n=1}^m\a_nx_n\biggr\V_\X\ge\biggl\V\sum_{n=1}^m\a_nz_n\biggr\V_\X.
$$
Condition 1) now follows since $\seq zn1$ is equivalent to $\seq {e^1}n1$.
\end{proof}

By Singer's result, any  $1$-dominated perturbation of $\seq {e^1}n1$ yields a copy of $\seq {e^1}n1$ consisting  of vectors that 
still possess large norm in $l_\infty$. Our next result, Theorem~2, is related to this observation. It states that any normalized unconditional 
basic sequence in $l_1$ consisting of vectors  with large $l_\infty$-norm is equivalent to $\seq{e^1}n1$. The proof of Theorem~2 is a 
consequence of the following fact, which, in turn, is a consequence of Szarek's refinement of the Khintchine inequality (c.f. \cite{LP}).

\begin{fact}\label{fc: abssum} 
The operator $\I\colon l_1\rightarrow l_2$ is absolutely summing; that is, $\I$ maps unconditionally convergent series into absolutely 
convergent series. Moreover, if the sequence $\seq xn1$ in $l_1$ has unconditional basis constant $C$, then for any sequence of scalars 
$\{\alpha_i\}_{1=1}^n$ we have
$$
  \biggl\V \sum_{i=1}^n \alpha_i x_i  \biggr\V_{l_1} \ge {1\over C\sqrt2}  \sum_{i=1}^n \V \alpha_i\I x_i \V_{l_2} 
$$
\end{fact}

\begin{theorem}\label{rm: bigl1seq} 
Let $\seq xn1$ be a normalized $K$-unconditional basic sequence in $l_1$. If, for some $1<p\le\infty$,
\begin{equation}\label{eq: www}
      \inf_{n\in\mathbb N} \V \II p x_n\V_{l_p}>0,
\end{equation}
then $\seq xn1$ is $(k,1)$-equivalent to the standard unit vector basis of $l_1$, where $k={1\over K\sqrt2}\inf_{n\in\mathbb N}\V\I x_n\V_{l_2}$.
\end{theorem}

\begin{proof}
First note that if $\seq xn1$  is a semi-normalized sequence in $l_1$ such that $\inf_{n\in\mathbb N} \V \II p x_n\V_{l_p}>0$ for some 
$1<p\le\infty$, then we necessarily have $\inf_{n\in\mathbb N} \V \II p x_n\V_{l_p}>0$ for each $1\le p\le\infty$. Indeed, this is obvious for $p=\infty$, 
while for $1<p<\infty$ and for $x=\sum_{n=1}^\infty \alpha_n e_n^1 \in l_1$ we have:
$$
  \V\II px\V_{l_p}^p=\sum_{i=1}^\infty|\alpha_i|^{p-1}|\alpha_i|\le \sum_{i=1}^\infty\V\II\infty x\V^{p-1}_{l_\infty}|\alpha_i|
   \le\V\II\infty x\V^{p-1}_{l_\infty}\cdot\V x\V_{l_1}. 
$$   
Thus, $\{\II\infty x_i\}_{i=1}^\infty$ is semi-normalized if both $\{\II px_i\}_{i=1}^\infty$ and $\{x_i\}_{i=1}^\infty$ are, whence the 
result follows.

In particular, equation~(2) holds for $p=2$. Using the triangle inequality and Fact~\ref{fc: abssum}, we have, for  any 
sequence of scalars $\{\alpha_i\}_{i=1}^n$, 
$$
  \sum_{i=1}^n |\alpha_i|\ge\biggl\V \sum_{i=1}^n\alpha_ix_i\biggr\V_{l_1}         
  \ge{1\over K\sqrt2}\sum_{i=1}^n\V\alpha_i\I x_i\V_{l_2}\ge k\sum_{i=1}^n|\alpha_i|,
$$
as desired.
\end{proof}
 
Finally, for a dual purpose, we present  Proposition~1. First, we see from Proposition~1, that the converse of Theorem~\ref{rm: bigl1seq} 
fails. Secondly, Proposition~1 sheds light on  Theorem~1. Note here that, for each $\epsilon>0$, any normalized unconditional 
basis of~$l_1$ is a  $(2+\epsilon)$-dominated perturbation of $\seq {e^1}n1$. It is thus natural to find the smallest
number $C$ such that the following statement holds: for each $\epsilon>0$, any normalized basis of~$l_1$ is a $(C+\epsilon)$-dominated 
perturbation of $\seq{e^1}n1$. From Proposition~1, we obtain that $C=2$.

\begin{proposition}\label{ex: eee} 
There is a semi-normalized unconditional basis of~$l_1$, hence a basis that is equivalent to the standard unit vector basis of $l_1$, satisfying
$$
  \inf_{n\in\mathbb N}\V \II\infty x_n\V_{l_\infty}=0,
$$  
where $\II\infty$ is the formal identity operator from $l_1$ into $l_\infty$.
\end{proposition}

To construct the required basis, we need the following lemma, whose proof, although certainly well-known, is provided for 
completeness.

\begin{lemma}\label{lm: jj}
Let $\{x_i\}_{i=1}^n$ be a basis of $l_1^n$ with the corresponding sequence of coefficient functionals $\{x_i^*\}_{i=1}^n$. Then 
$\{x_i\}_{i=1}^n$ is $(k_1,k_2)$-equivalent to the standard unit vector basis of $l_1^n$, where 
$k_1^{-1}=\max_{1\le i\le n}\sum_{j=1}^n|x_j^*(i)|$ and $k_2=\max_{1\le j\le n}\V x_j\V_{l_1^n}$.
\end{lemma}

\begin{proof}
We need to show that 
$$
  k_1\sum_{i=1}^n|\a_i|\le\biggl\V\sum_{i=1}^n\a_ix_i\biggr\V_{l_1^n}\le k_2\sum_{i=1}^n|\a_i|.
$$
for any sequence of scalars $\{\alpha_i \}_{i=1}^n$. The right hand inequality follows from the triangle inequality. Towards proving 
the left hand inequality, let $\{ \alpha_i \}_{i=1}^n$ be given and write
$$
  z=\sum_{i=1}^n \alpha_i x_i=\sum_{i=1}^n \beta_i e_i,
$$ 
where $\{e_i\}_{i=1}^n$ denotes the standard unit vector basis of $l_1^n$. Then for each $j=1,2,\ldots,n$, we have 
$\alpha_j=\sum_{i=1}^n\beta_i x_j^*(i)$; whence,
$$
  \sum_{j=1}^n |\alpha_j|=\sum_{j=1}^n\biggl|\sum_{i=1}^n\beta_ix_j^*(i)\biggr|  
  \le\sum_{i=1}^n|\beta_i|\sum_{j=1}^n|x_j^*(i)|\le{1\over k_1}\sum_{i=1}^n|\beta_i| ={1\over k_1}\V z\V_{l_1^n},
$$
as desired.   
\end{proof}

\begin{proof}[Proof of Proposition~1]
To construct the sequence heralded by Proposition~\ref{ex: eee}, it suffices to find, for each integer $n>2$, a basis 
$\{x^n_i\}_{i=1}^n$ of $l_1^n$ such that both of the following statements hold:
\begin{itemize}
      \item[1)]$\{x^n_i\}_{i=1}^n$ is $(\ts{1\over5},2)$-equivalent to the standard unit vector basis of~$l_1^n$.
      \item[2)]$\max\limits_{1\le i\le n} |x_1(i)|  =1/n$.
\end{itemize}
Towards this end, let $n\in \mathbb N$ and let $\{e_i\}_{i=1}^n$ denote the standard unit vector basis of $l_1^n$. Define
\begin{eqnarray*}
       x_1&=&\sum_{i=1}^n{\ts {1\over n}} e_i                                                 \\
       x_i&=&e_1+e_i \quad{ \rm for\ } 2\le i\le n                                            \\
\end{eqnarray*}
(we drop the superscripts for notational clarity). Then we have condition 2), and $\{x_i\}_{i=1}^n$ is a linearly independent sequence 
with coefficient functionals $\{x^*_i\}_{i=1}^n$ given by:
\begin{eqnarray*}
      x^*_1&=&{\ts{-n\over n-2}}e_1+\sum_{i=2}^n{\ts{n\over n-2}}e_i                          \\
      x^*_j&=&{\ts{1\over n-2}}e_1+{\ts{n-3\over n-2}}e_j +\sum_{i\notin\{j,1\}}{\ts{-1\over n-2}}e_i\quad{\rm for\ }2\le j\le n.  
\end{eqnarray*}
Using Lemma~\ref{lm: jj}, it is easily verified that $\{x_i\}_{i=1}^n$ is $(1/5,2)$-equivalent to $\{e_i\}_{i=1}^n$.
\end{proof}

\noindent{\bf Acknowledgements}
The author would like to thank the members of the Functional Analysis Seminar at the University of South 
Carolina for their interest in this work, and for their helpful suggestions towards simplifying and generalizing the results of this paper.


\enddocument